\documentclass{article}

\usepackage{times}
\usepackage[T1]{fontenc}
\usepackage[utf8]{inputenc}
\usepackage{array}
\usepackage{multirow}
\usepackage{float}
\usepackage{url}
\usepackage{amsthm}
\usepackage{amsmath}
\usepackage{amssymb}
\usepackage{algorithm}
\usepackage{algorithmic}
\usepackage{mathtools}
\usepackage{graphicx}
\usepackage{color}
\usepackage[english]{babel}
\usepackage{dsfont} 

\newcommand\Mark[1]{\textsuperscript#1}

\newcommand{\G}{\mathcal{G}} 
\newcommand{\E}{\mathcal{E}} 
 
\newcommand{\Lap}{\mathcal{L}} 
\newcommand{\V}{\mathcal{V}} 
\newcommand{\Om}{\mathcal{O}}

\newcommand{\x}{u}

\newcommand{\Rbb}{\mathbb{R}}

\newtheorem{theorem}{Theorem}

\providecommand{\keywords}[1]{\textbf{\textit{Index terms---}} #1}

\begin{document}


\begingroup
\centering
{\LARGE Accelerated filtering on graphs using Lanczos method\\[1.5em]
\large Ana {\v S}u{\v s}njara\Mark{1}, Nathanaël Perraudin\Mark{2}, Daniel Kressner\Mark{1}, and
Pierre Vandergheynst\Mark{2}}\\[1em]
\begin{tabular}{*{2}{>{\centering}p{.25\textwidth}}}
 \Mark{1} ANCHP & \Mark{2}LTS2 \tabularnewline
Ecole Polytechnique Fédérale de Lausanne
EPFL & Ecole Polytechnique Fédérale de Lausanne
EPFL  \tabularnewline
\end{tabular}\par
\endgroup

\begin{abstract}
Signal-processing on graphs has developed into a very active field of research during the last decade.
In particular, the number of applications using frames constructed from graphs, like wavelets on graphs, has substantially increased. To attain scalability for large graphs, fast graph-signal filtering techniques are needed. 
In this contribution, we propose an accelerated algorithm based on the Lanczos method that adapts to the Laplacian spectrum without explicitly computing it. The result is an accurate, robust, scalable and efficient algorithm. 
Compared to existing methods based on Chebyshev polynomials, our solution achieves higher accuracy without increasing the overall complexity significantly. Furthermore, it is particularly well suited for graphs with large spectral gaps.
\footnote{The work of A. {\v S}u{\v s}njara has been supported by the SNSF research project \textit{Low-rank updates of matrix functions and fast eigenvalue solvers}. The work of Nathanaël Perraudin has been supported by the SNF research project \textit{Towards Signal Processing on Graphs}, grant number: 2000\_21/154350/1 .}

\end{abstract}

\section{Introduction}
\begin{figure}[htb!]
\begin{center}
\includegraphics[width=0.9\textwidth]{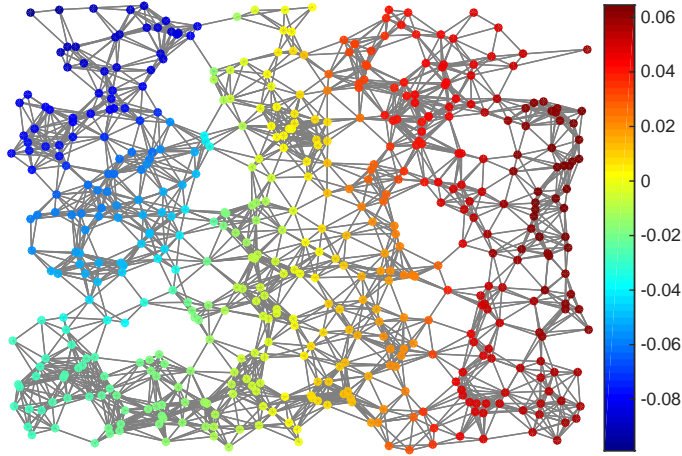} 
\end{center}
\caption{Synthesized sensor network. Colored circles represent vertices and gray links represent edges. Values of the signal are displayed with different colors.}
\label{fig:sensor_graph}
\vspace{-10pt}
\end{figure}

Graphs conveniently represent data on irregular geometric structures as they arise in numerous application domains such as social, energy, transportation or neuronal networks. In all of these fields, different pieces of information are connected to each other and these connections can be modelled by a graph. For every link, an edge is drawn with an associated weight that represents the similarity between the two elements (vertices) it connects. For example, in a sensor network acquiring environmental measurements, it makes sense to choose the weight inversely proportional to the physical distance. 
A data signal lives on nodes and can be visualized as a collection of samples. We refer to those as graph-signals, see Figure~\ref{fig:sensor_graph} for an example.

This framework, called graph signal-processing, has recently become a general field of research \cite{shuman2013emerging,shuman2013vertex} and applications of graph signal-processing can be found in many different areas. In machine vision, automatic text classification or more generally machine learning, graphs are used to represent similarities between data points and result in algorithms for semi-supervised learning \cite{zhu2003semi,smola2003kernels,belkin2004regularization,zhou2004learning}. Graph signal-processing is often used to leverage intrinsic links when processing data. For instance, many applications in image processing benefit from the use of graphs that describe connections between local patches in the image (e.g.~\cite{peyre2008non,zhang2008graph} and references therein). 

Spectral graph filtering is one of the basic building blocks in the applications discussed above. 
Classical filtering techniques rely on the Fourier transform, which can be done inexpensively with a cost of $\Om(N\log(N))$.  Extending this transform to graphs requires the diagonalisation of the graph Laplacian, which becomes
prohibitively expensive for larger graphs. To overcome this issue, an efficient approximate method to perform spectral graph filtering was proposed in~\cite{hammond2011wavelets}. This method, based on Chebyshev polynomial approximation, only requires multiple application of the Laplacian operator and thus leads to a scalable, efficient memory-saving and error-controlled algorithm.

Instead of Chebyshev polynomial approximation, we propose to use the Lanczos method to perform spectral graph filtering. The resulting algorithm is also scalable, in addition to being superior in accuracy to its predecessor.

The rest of this article is organized as follows. In Section 2, we summarize the basics of signal-processing on graphs, including the definition of graph filtering. We also recall the method of Hammond et al.~\cite{hammond2011wavelets}. Section 3 describes the Lanczos method and its application to graph filtering. Numerical experiments are presented in Section 4.

\section{Signal-processing on graphs}

\subsection{Graph nomenclature}

We consider a weighted undirected graph $\G=\{ \V,\E,\mathcal{W}\}$ with a set of vertices $\V$, a set of edges $\E$,
and a weight function $\mathcal{W} : \V \times \V \rightarrow \Rbb$. The vertices are indexed from $1,\dots, N=|\V|$ and each entry  of the weight matrix $W \in \mathbb R^{N\times N}$ contains the weight of the edge connecting the corresponding vertices: $W_{i,j} = \mathcal{W}(v_i,v_j)$. If there is no edge between two vertices, the weight is set to $0$. It is assumed that $W_{i,j} = W_{j,i}$, that is, $W$ is a symmetric matrix. For a vertex $v_i\in \V$, the degree $d(i)$ is defined as the sum of the weights of incident edges: $d(i)=\sum_{j=1}^N W_{i,j}$. 

In this framework, a graph signal is defined as a function $s: \V \rightarrow  \mathbb{R}$ assigning a value to each vertex. It is convenient to consider a signal $s$ as a vector of size $N$ with the $i^{\mathrm{th}}$ component representing the signal value at the $i^{\mathrm{th}}$ vertex. 

One of the most fundamental concepts for weighted undirected graphs is the (combinatorial) graph Laplacian $\Lap$ 
defined as
$\Lap = D-W,$ where $D$ is the diagonal degree matrix with diagonal entries $D_{ii}=d(i)$.
Alternative definitions of graph Laplacians include the normalized Laplacian $\Lap_n = D^{\frac{1}{2}} \Lap D^{\frac{1}{2}} = D^{\frac{1}{2}} (D-W) D^{\frac{1}{2}}$.
For simplicity, we restrict ourselves to $\Lap$, but the method presented in this paper easily extends to other Laplacians. 


The Laplacian $\Lap$ is always symmetric positive semi-definite and can thus be decomposed as 
$$\Lap = U\Lambda U^*,$$
where $U = [\x_0, \ldots, \x_{N-1}]$ is an orthogonal matrix and $U^*$ denotes its transpose.
Without loss of generality, we order the set of eigenvalues as follows: $0=\lambda_0 \leq \lambda_1 \leq ... \leq \lambda_{N-1} = \lambda_{\rm max}$;
see~\cite{chung1997spectral} for more details on spectral graph theory. 
The matrix $U$ defines the graph Fourier basis~\cite{shuman2013emerging,shuman2013vertex}, leading to the graph Fourier transform $\hat{s} = U^*s$ and its inverse $s = U\hat{s}$.


\subsection{Graph filters}
In the classic setting, applying a filter to a signal is performed by convolution, which corresponds to point-wise multiplication in the spectral domain. 
Similarly, filtering a graph-signal is peformed by multiplication with a filter in the graph Fourier domain. A graph filter is defined by a continuous function $g: \Rbb_+ \rightarrow \Rbb$. To obtain its discrete coefficients, this function is evaluated at each eigenvalue: $g(\lambda_\ell)$ for $\ell  = 0,\ldots,N-1$. The filtering operation then corresponds to $\hat{s}'(\ell) =  g(\lambda_\ell) \cdot \hat{s}(\ell)$, where $s'$ is the filtered signal. Equivalently, using matrix notation, we have
\begin{equation} \label{eq:filteredsignal}
s' = U g(\Lambda) U^* s\text{,}
\end{equation}
where $g(\Lambda)$ is the diagonal matrix containing the coefficients $g(\lambda_\ell)$ on the diagonal. 
In terms of matrix functions~\cite{Higham2008}, the relation~\eqref{eq:filteredsignal} can be compactly expressed as
$s' = g(\Lap) s$ with $g(\Lap) := U g(\Lambda) U^{*}$.


\subsection{Fast filtering via Chebyshev polynomials} \label{sec:fastfiltering}

The graph filtering operation described above is based on the graph Fourier transform. Unfortunately, the graph Fourier basis needed for performing this transform requires the diagonalization of the graph Laplacian, which takes $\Om(N^3)$ operations and $\Om(N^2)$ memory when using standard techniques. This is feasible for graphs with only a few thousand vertices. To be able to tackle problems of larger size, more efficient methods are needed, one of which is presented in \cite{hammond2011wavelets}
and summarized in this section.

\textit{Filtering in the vertex domain. } To avoid the Fourier transform, we perform the filtering operation in the vertex domain using only the Laplacian operator. Applying this operator corresponds to multiplying the signal in the spectral domain with the eigenvalues:
$$
\widehat{\Lap s} = \Lambda \hat{s}\text{.}
$$
This is equivalent to filtering with $g(x) = x$. Using this relation recursively and exploiting linearity, we can apply any polynomial filter $g(x) = a_0 + a_1 x  +\dots + a_K x^M$ to a signal $s$ with the following formula:
\begin{equation} \label{eq:polynomial_filter}
s' = U g(\Lambda) U^* s = \left(a_0 I + a_1 \Lap + \dots + a_M \Lap^M \right)s.
\end{equation}

\textit{Chebyshev polynomial approximation. }
The ability to apply polynomial filters efficiently suggests to 
approximate a given filter function with a suitable polynomial.
For approximating functions on real intervals, Chebyshev polynomials are usually the preferred choice because of numerical stability considerations and the fact that they can be evaluated efficiently by three-term recurrences.
We refer to~\cite{hammond2011wavelets} for a more detailed discussion on the choice of Chebyshev polynomials in signal-processing on graphs and to, e.g.,~\cite{phillips2003interpolation} for an introduction to polynomial approximation.

The $m$th Chebyshev polynomial $T_m(y)$ is generated using the recurrence relation $T_m(y) = 2x T_{m-1}(y) - T_{m-1}(y)$ with $T_0(y) = 1$ and $T_1(y) = y$. For $y\in[-1,1]$, these polynomials possess the following well-known properties:
\begin{enumerate}
  \item they admit the closed form expression $T_m(y) = \cos (m \arccos (y))$;
  \item they are bounded, i.e., $T_m(y) \in [-1,1]$; 
  \item they form an orthogonal basis of $L^2\left([-1,1], \frac{\mathrm{d}y}{\sqrt{1-y^2}}\right)$.
\end{enumerate}
The third property implies that every function $h\in L^2\left([-1,1], \frac{\mathrm{d}y}{\sqrt{1-y^2}}\right)$ admits a convergent Chebyshev series
\begin{equation*}
h(y) = \frac{1}{2} c_0 + \sum_{m=1}^\infty c_m T_m(y),
\end{equation*}
with the Chebyshev coefficients
\begin{equation*}
c_k  
= \frac{2}{\pi}\int_{-1}^1 \frac{T_m(y)h(y)}{\sqrt{1-y^2}}\mathrm{d}y 
= \frac{2}{\pi}\int_{0}^\pi \cos (k\theta) h(\cos (\theta)) \mathrm{d}\theta.
\end{equation*}
Since our filter $g$ is evaluated on the eigenvalues of the Laplacian, we need to map the interval $[-1,1]$ to the interval
$[0,\lambda_{\max}]$ using the transformation $x = \frac{\lambda_{\max}}{2}(y+1)$. Defining $\tilde{T}_m(x) = T_m\left( \frac{2x}{\lambda_{\max}}-1 \right)$ we obtain
\begin{equation} \label{eq:cheb_approx}
g(x) = \frac{1}{2} c_0 + \sum_{m=1}^\infty c_m \tilde{T}_m(x)
\end{equation}
for $x\in [0,\lambda_{\max}]$, with 
\begin{equation*}
c_m  
= \frac{2}{\pi}\int_{0}^\pi \cos (m\theta) g\left(\frac{\lambda_{\max}}{2}\left(\cos (\theta)+1 \right) \right) \mathrm{d}\theta.
\end{equation*}

\textit{Fast filtering algorithm. }
At this point, we can derive the iterative algorithm for filtering a signal $s$ with $g$. The recurrence relation for the transformed Chebyshev polynomials becomes 
$\tilde{T}_m(x)
= 2\left( \frac{2x}{\lambda_{\max}}-1 \right) \tilde{T}_{m-1}(x)
- \tilde{T}_{m-2}(x)$. On the matrix level, this yields, using \eqref{eq:polynomial_filter}:
\begin{equation*}
\tilde{T}_m(\Lap) s=  2\left( \frac{2\Lap}{\lambda_{\max}}-I \right) \tilde{T}_{m-1}(\Lap)s
- \tilde{T}_{m-2}(\Lap)s.
\end{equation*}
Combined with \eqref{eq:cheb_approx}, this finally leads to the following expression for filtering a signal $s$:
\begin{equation*}
s' = g(\Lap) = \frac{1}{2} c_0 I s + \sum_{m=1}^\infty c_m \tilde{T}_m(\Lap)s.
\end{equation*}
When implemented, we truncate this sum at a defined order $M$. Assuming that $|\E|>N$, the computational cost of this algorithm scales linearly with the number of edges $\Om(M|\E|)$. In most applications, the Laplacian is sparse,
$|\E|\ll N^2$, which results in a fast algorithm. Moreover, apart from storing the Laplacian, the additional memory consumed by this algorithm is only $4 N$.

\section{Accelerated filtering using Lanczos} \label{sec:lanczos}

Given the graph Laplacian $\Lap\in\mathbb{R}^{N\times N}$ and a nonzero vector $s\in\mathbb{R}^N$, the Lanczos method~\cite{golub2013matrix} shown in Algorithm~\ref{alg:alg1} below computes an orthonormal basis $V_M = [v_1,\ldots,v_M]$ of the Krylov subspace $K_M(\Lap,s)  = \operatorname{span}\{s, \Lap s,\ldots,\Lap^{M-1} s \}$. 
The computational cost of Algorithm is $\Om(M\cdot |\E|)$.
The storage of the basis $V_M$ requires $M N$ additional memory, which can be avoided using two passes of the algorithm or
restart techniques; see, e.g.,~\cite{Frommer2014} for more details.
   \begin{algorithm}
    \caption{\text{Lanczos method}}
    \label{alg:alg1}
    \renewcommand{\algorithmicrequire}{\textbf{Input:}}
    \renewcommand{\algorithmicensure}{\textbf{Output:}}

    \begin{algorithmic}[1]
   
    \REQUIRE Symmetric matrix $\Lap \in \mathbb{R}^{N\times N}$, vector $s\neq 0$, $M\in \mathbb{N}$.
    \ENSURE $V_M = [v_1,\ldots,v_M]$ with orthonormal columns, scalars $\alpha_1,\ldots,\alpha_M \in \mathbb R$ and $\beta_2,\ldots,\beta_M \in \mathbb R$.
 
   \STATE $v_1 \gets s/ \Vert s\Vert_2$
   \FOR{$j = 1,2,\ldots,M$} 
        \STATE {$w = \Lap v_j$}
        \STATE {$\alpha_j = v_j^{*}w$}
        \STATE {$\tilde{v}_{j+1} = w - v_j\alpha_j$}
          \IF {$j > 1$}
          \STATE {$\tilde{v}_{j+1} \gets \tilde{v}_{j+1} - v_{j-1}\beta_{j-1}$}
          \ENDIF
       \STATE {$\beta_{j} = \Vert \tilde{v}_{j+1}\Vert_2$}
       \STATE $v_{j+1} = \tilde{v}_{j+1}/\beta_{j}$     
   \ENDFOR     
   \end{algorithmic}
   \end{algorithm}

The scalars produced by Algorithm~\ref{alg:alg1} can be arranged into 
a symmetric tridiagonal matrix $H_M \in \mathbb R^{M\times M}$ satisfying 
\begin{equation*}
\label{eqref:projection}
V_M^{*} \Lap V_M = H_M = \begin{bmatrix}
\alpha_1 &\beta_2\\
\beta_2 &\alpha_2 &\beta_3\\
 &\beta_3 &\alpha_3 &\ddots &\\
& &\ddots & \ddots & \beta_{M} \\
& & &\beta_{M} &\alpha_{M} 
\end{bmatrix}\text{.}
\end{equation*}
In floating-point arithmetics, the orthogonality of the basis produced by Algorithm~\ref{alg:alg1} may get quickly lost and reorthogonalization is needed~\cite{cullum2002lanczos}.

Given a continuous function $g: [0,\lambda_{\max}] \to \mathbb R$ and a vector $s$, the following approximation
to $g(\Lap)s$ was proposed by Gallopoulos and Saad in~\cite{gallopoulos1992efficient}:
\begin{equation}
\label{eq:final_step}
g(\Lap)s \approx \Vert s \Vert_2 V_M g(H_M)e_1\coloneqq g_M\text{,}
\end{equation}
where $e_1 \in \Rbb^M$ is the first unit vector. Because of eigenvalue interlacing,
the eigenvalues of $H_M$ are contained in  the interval $[0,\lambda_{\max}]$ and
hence the expression $g(H_M)$ is well-defined. Typically, 
$M \ll N$, rendering the evaluation of $g(H_M)$ inexpensive. The overall cost of
our Lanczos-based approximation of graph-signal filtering, which consists of applying Algorithm~\ref{alg:alg1}
and evaluating~\eqref{eq:final_step}, is therefore between 
$\Om (M \cdot |\E|)$ and $\Om (M \cdot |\E| + M^2N)$, depending on how the reorthogonalization is performed.

\textit{Approximation quality.} The following theorem provides some insight into the accuracy of
the approximation $g_M$ obtained by the Lanczos method.
\begin{theorem}[{\cite[Corollary 3.4]{Guttel2013a}}]
\label{thm1}
Let $\Lap\in\mathbb{R}^{N\times N}$ be symmetric with eigenvalues contained in the interval $[ 0, \lambda_{\max}]$ and let $g:[ 0, \lambda_{\max}]\rightarrow \mathbb R$ be continuous.
Then 
\begin{align*}
&\Vert g(\Lap)s - g_M\Vert_2 
\leq 2\Vert s \Vert_2 \cdot  \underset{p\in \mathcal{P}_{M-1}}{\operatorname{min}}\, \underset{z \in \left[0, \lambda_{\max} \right]}{\operatorname{max}} \vert g(z) - p(z)\vert\text{,}
\end{align*}
where $\mathcal{P}_{M-1}$ denotes all polynomials of degree at most $M-1$.
\end{theorem}
Theorem~\ref{thm1} shows that the error is bounded by the best polynomial approximation~\cite{phillips2003interpolation} of $g$ on $[ 0, \lambda_{\max}]$.
Compared to the Chebyshev approximation from Section~\ref{sec:fastfiltering}, this shows that -- up to a multiple of two --
the Lanczos-based approximation $g_M$ can be expected to provide at least the same accuracy.

However, the Lanczos-based approximation can sometimes be expected to perform significantly better because of its ability to adapt to the
eigenvalues of $\Lap$. This phenomenon is well-understood for Krylov subspace approximations to solutions of linear systems (see, e.g.,~\cite[Section 3.1]{Greenbaum1997})
and extends to matrix functions as well. To illustrate such a situation, suppose that the eigenvalue $0$ of $\Lap$ is well separated from the other eigenvalues
contained in the (smaller) interval $[\lambda_1,\lambda_{\max}]$.
Then the eigenvalues of $H_M$ can be expected to exhibit a similar behavior~\cite[Chapter 6]{Saad2011}.

By choosing a polynomial that interpolates the eigenvalue $0$ exactly and approximates $g$ on $[\lambda_1,\lambda_{\max}]$, the bound of Theorem~\ref{thm1} can approximately be replaced by
\[
 \Vert g(\Lap)s - g_M\Vert_2 \lesssim 2\Vert s \Vert_2 \cdot  \underset{p\in \mathcal{P}_{M-2}}{\operatorname{min}}\, \underset{z \in \left[\lambda_1, \lambda_{\max} \right]}{\operatorname{max}} \vert g(z) - p(z)\vert\text{.}
\]
At the expense of losing a degree of freedom in the choice of the polynomial, the width of the approximation interval becomes smaller,
which in turn leads to significantly improved approximation rates~\cite{phillips2003interpolation}. In contrast, the Chebyshev approximation cannot
adapt to the eigenvalues of $\Lap$ and hence its approximation rate is driven by the larger interval $[0,\lambda_{\max}]$.

\textit{Stopping criteria.} Ideally, $M$ should be chosen to be the smallest integer such that 
\begin{equation*}
\Vert e_M \Vert_2 = \Vert g(\Lap)s - g_M \Vert_2\leq \epsilon 
\end{equation*} 
holds for some prescribed accuracy $\epsilon>0$. 
As proposed in~\cite{vandeneshof2006preconditioning}, we estimate $e_M$ as the difference of two (consecutive) approximations:
\begin{equation}
\label{eq:err_est}
e_M \approx g_{M+j} - g_{M}
\end{equation}
for a small value of $j$.

\section{Numerical experiments}

In this section we present numerical experiments to demonstrate the numerical behaviour of the Lanczos method and compare it with the Chebyshev method. All experiments were performed with the signal processing toolbox GSPBox~\cite{perraudin2014gspbox} and can be downloaded at \url{https://lts2.epfl.ch/lanczos-filtering/}. 

\textit{Example 1.} We first consider the sensor graph and the Erd{\H o}s--R{\'e}nyi random graph~\cite{erdHos1959random}, with $N = 500$ nodes. In the latter case, each edge is included in the graph with probability $p = 0.04$. As a filterbank, we use a collection of translated windows $g(t) = \sin(0.5\pi\cos(\pi t)^2)\chi_{[-1/2,1/2]}$, where $\chi_{I}$ is the characteristic function of $I$ \cite{wesfreid1993adapted}. 
We choose to adapt the filterbank to the graph spectrum, since non-adapted filters lead to very coherent frames for graphs with eigenvalues not uniformly spread along the spectrum~\cite{shuman2013spectrum}. 

\begin{figure}[htb!]
\begin{center}
\includegraphics[width=0.45\textwidth]{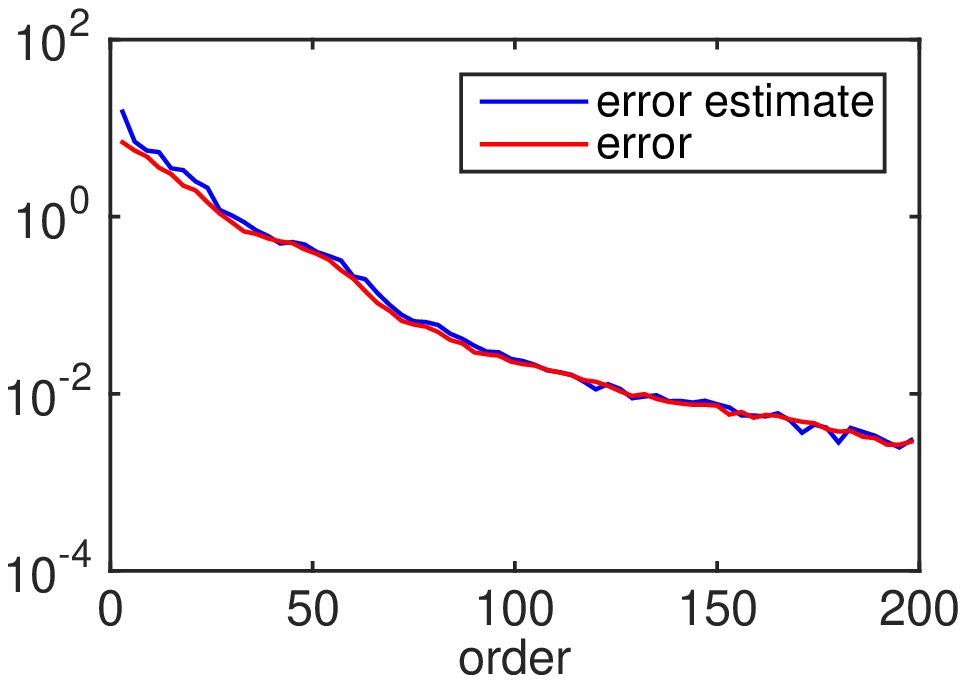} 
\includegraphics[width=0.45\textwidth]{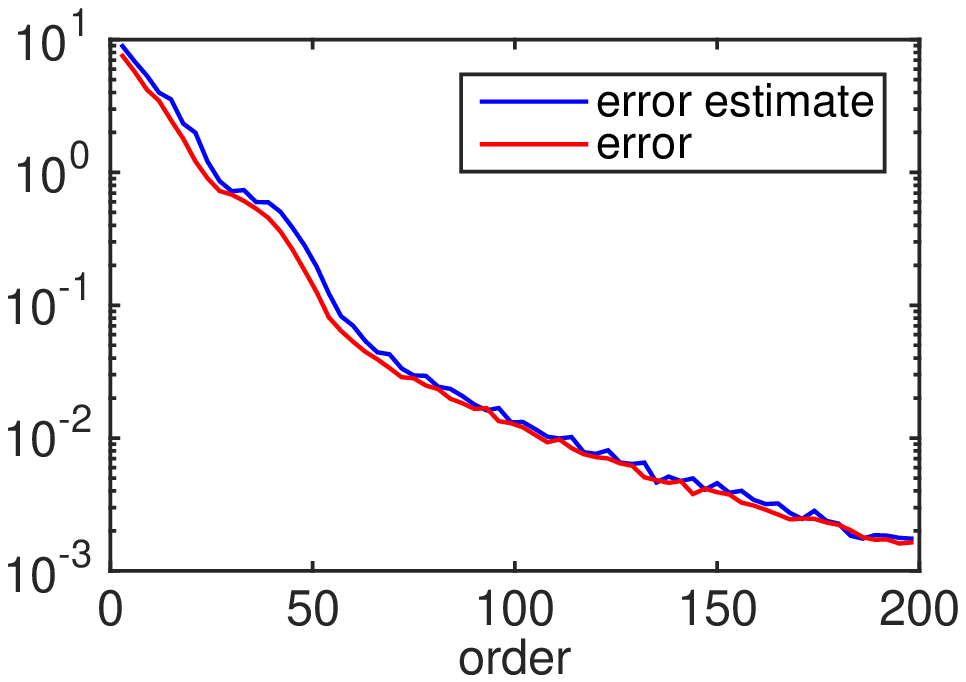} 
\caption{Sensor graph (left) and Erd{\H o}s – R{\'e}nyi graph (right). Comparison of the approximation error in the Lanczos method and the estimate~\eqref{eq:err_est}.}
\label{fig:error_vs_est}
\end{center}
\begin{center}
\includegraphics[width=0.45\textwidth]{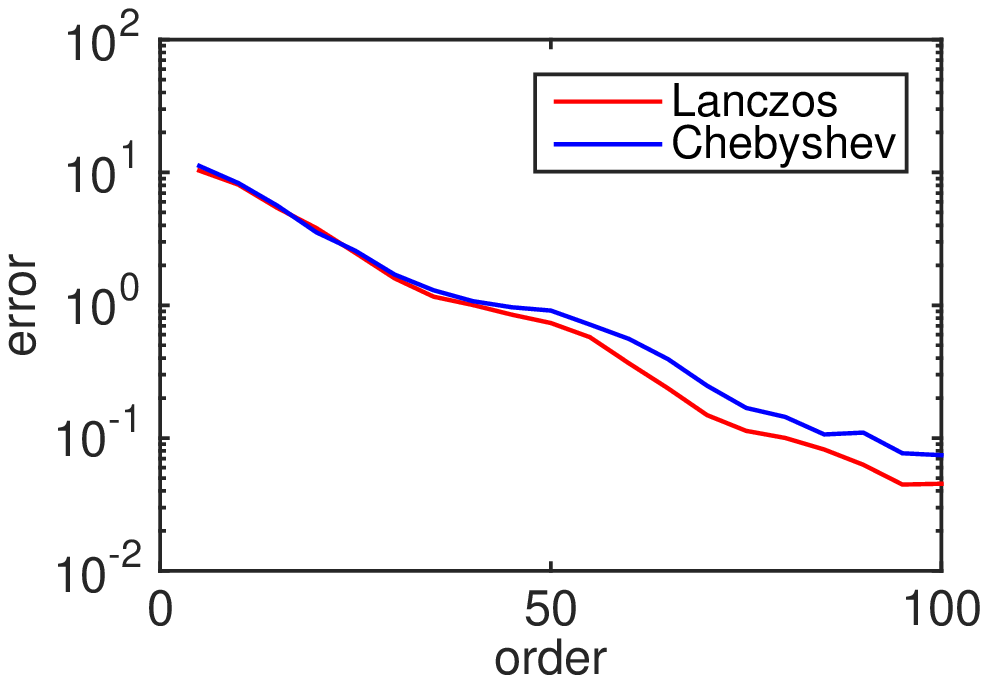} 
\includegraphics[width=0.45\textwidth]{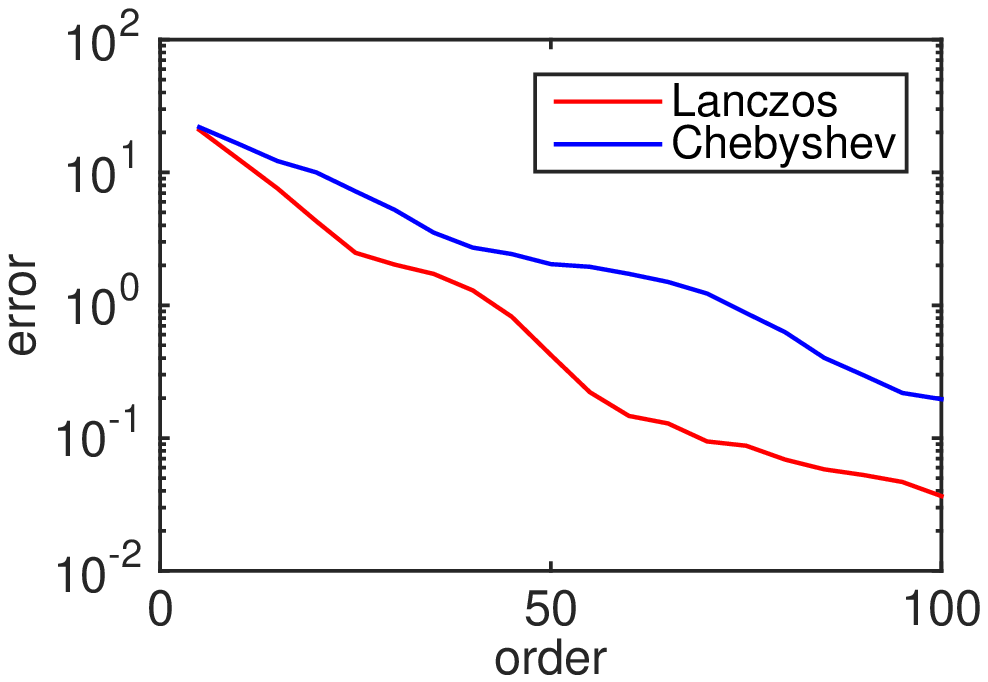} 
\caption{Sensor graph (left) and Erd{\H o}s – R{\'e}nyi graph (right). Comparison of the approximation error using the Chebyshev method and the Lanczos method with respect to the order of approximation.}
\label{fig:error_vs_order}
\end{center}
\vspace{-10pt}
\end{figure}

First, we test the reliability of the error estimate~\eqref{eq:err_est} for $j = 3$. Figure~\ref{fig:error_vs_est} confirms that the estimate
$\|g_{M+3} - g_{M}\|_2$ closely follows the true error $\|e_M\|_2$.

\begin{figure}[ht!]
\begin{center}
\includegraphics[width=0.45\textwidth]{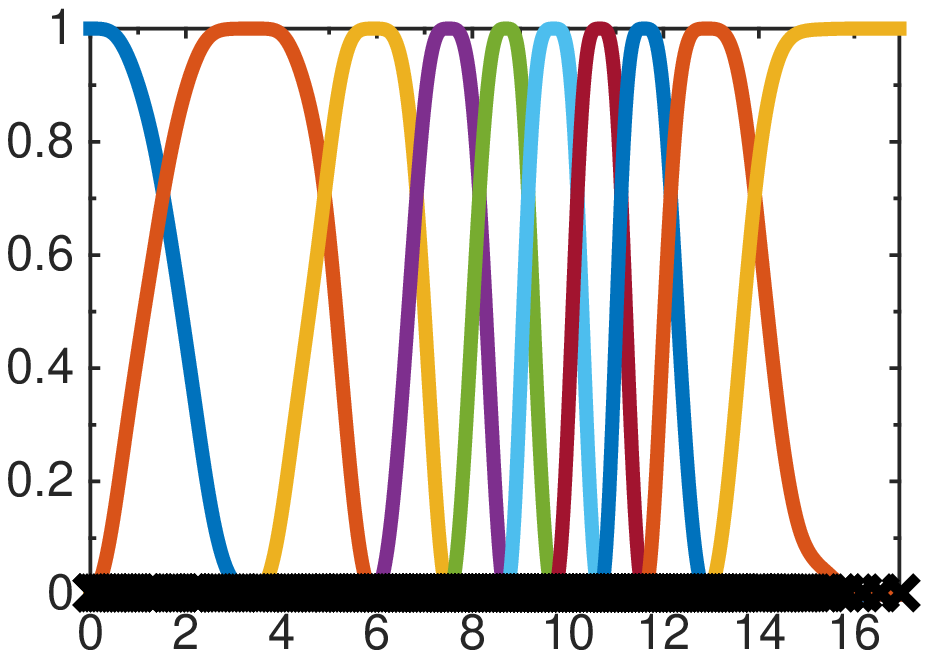} 
\includegraphics[width=0.45\textwidth]{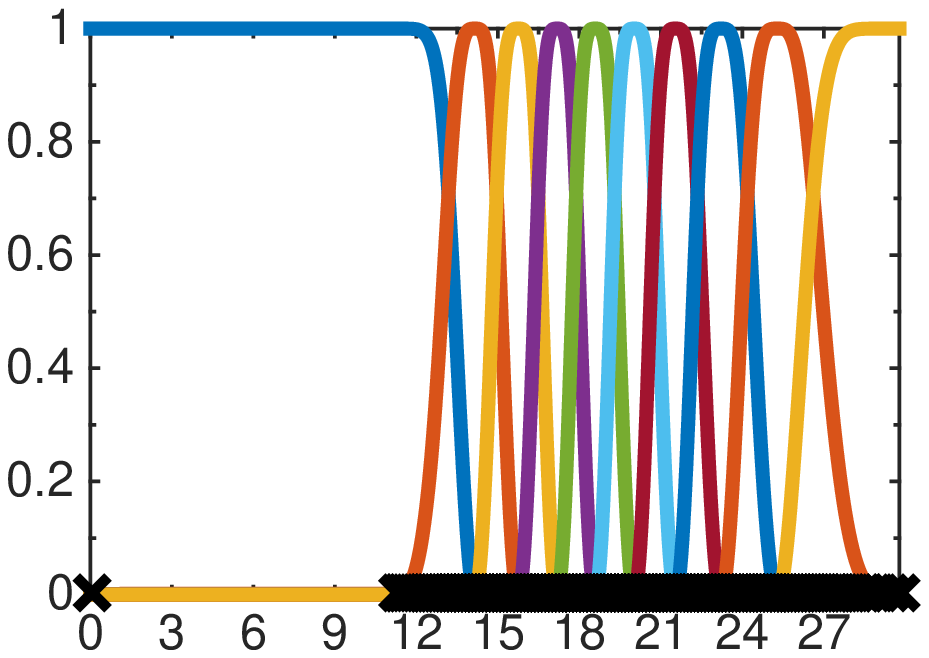} 
\caption{Sensor graph (left) and Erd{\H o}s – R{\'e}nyi graph (right). Filterbanks adapted to spectra of graph Laplacians.}
\end{center}
\begin{center}
\includegraphics[width=0.45\textwidth]{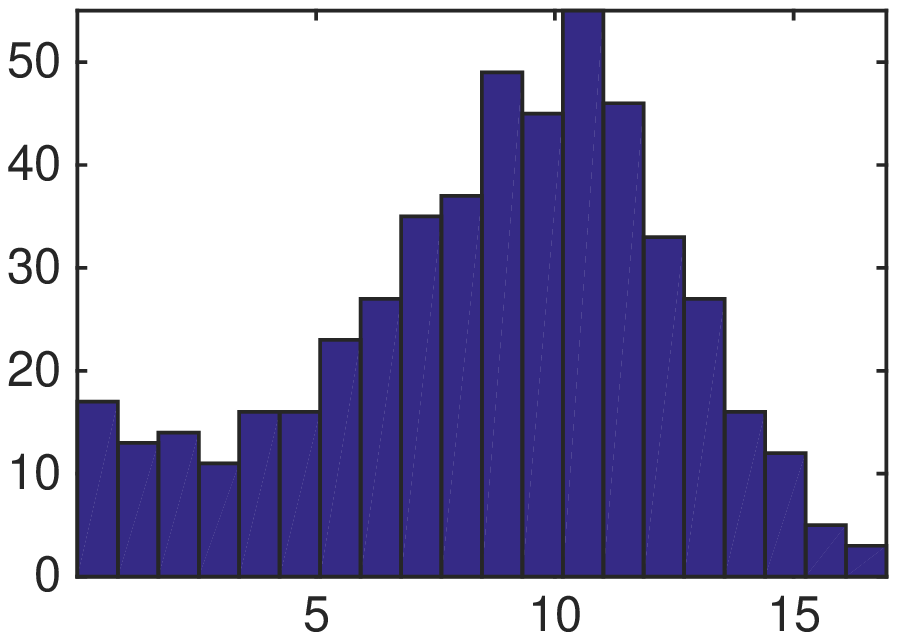} 
\includegraphics[width=0.45\textwidth]{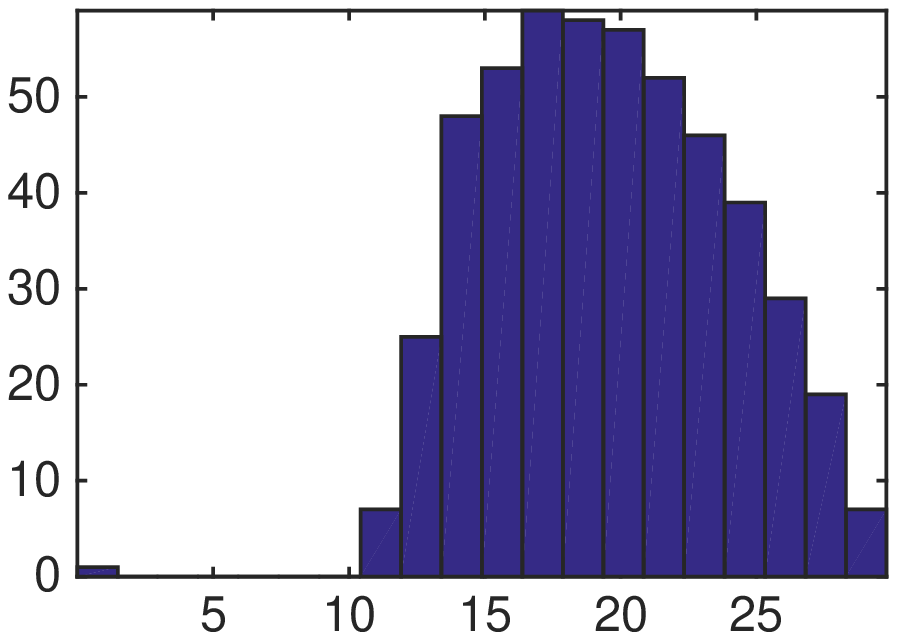} 
\caption{Sensor graph (left) and Erd{\H o}s – R{\'e}nyi graph (right). Distribution of graph Laplacians' eigenvalues.}
\label{fig:eig_hist}
\end{center}
\vspace{-10pt}
\end{figure}
In Figures~\ref{fig:error_vs_order} and~\ref{fig:eig_hist} we observe that the Lanczos method yields better approximation, especially in the case when the spectrum of graph Laplacian is not uniformly distributed and has a large (relative) spectral gap, as predicted in Section~\ref{sec:lanczos}.     

\textit{Example 2.} The purpose of this example is to investigate how the approximation by the Chebyshev polynomial method and the Lanczos method depends on the edge sparsity of the Erd{\H o}s – R{\'e}nyi graph. We consider graphs with $N = 1000$ nodes. The order of Chebyshev polynomial and Lanczos method are set to $M = 30$. As a filterbank, we the use mexican hat wavelet, with a mother window 
\begin{equation*}
g_h(\lambda_\ell) =  \lambda_\ell \cdot \operatorname{exp}(-\lambda_\ell^2)\text{,}
\end{equation*}
with $\lambda_\ell$ eigenvalues of graph Laplacian. In order to obtain a complete filterbank, we add a low pass filter
$g_l(\lambda_\ell) =  \operatorname{exp}(-\lambda_\ell^4)$.
   
\begin{figure}[ht!]
\begin{center}
\includegraphics[width=0.45\textwidth]{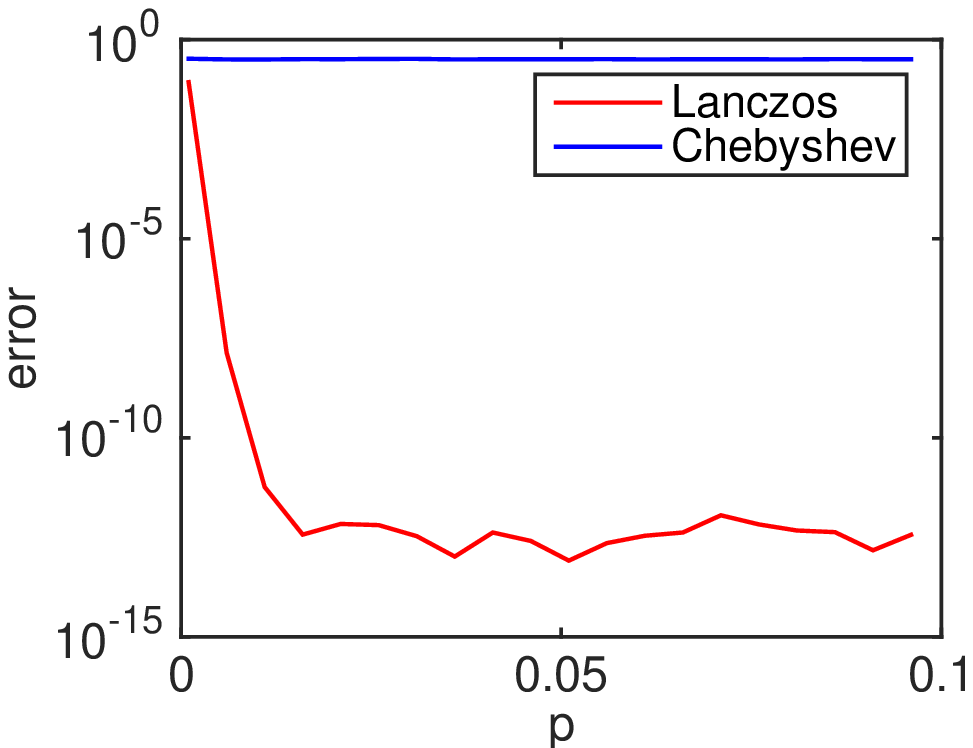}
\includegraphics[width=0.45\textwidth]{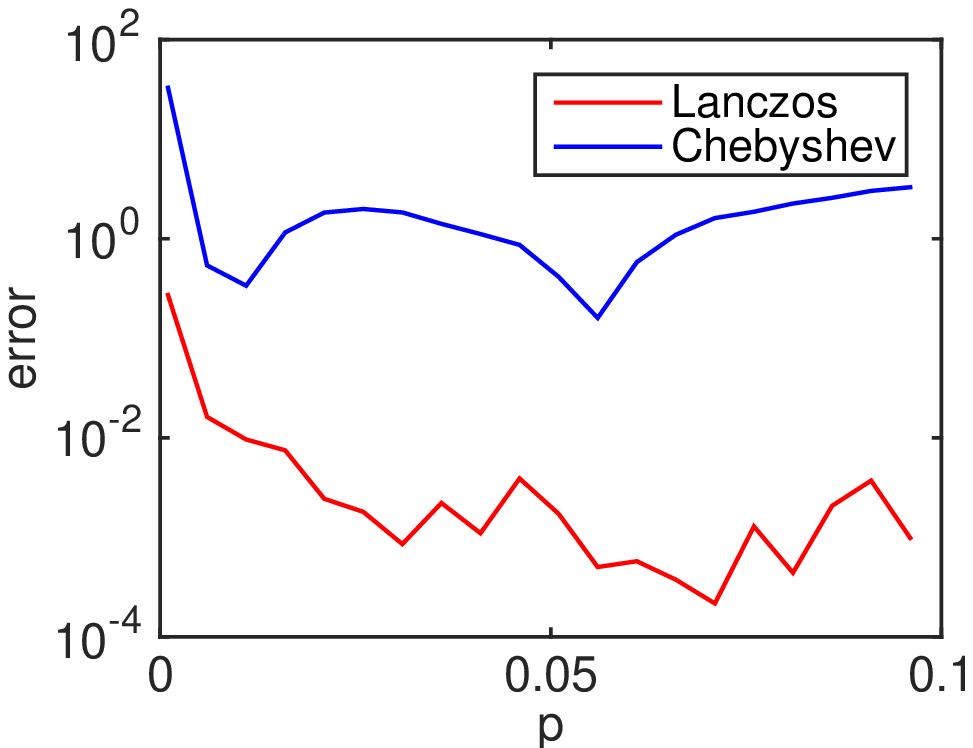} 
\caption{Comparison of the approximation error using the Chebyshev method and the Lanczos method with respect to $p$, with non-adapted filterbank (left) and adapted filterbank (right).}
\label{fig:error_vs_p}
\end{center}
\vspace{-10pt}
\end{figure}
\begin{figure}[ht!]
\begin{center}
\includegraphics[width=0.45\textwidth]{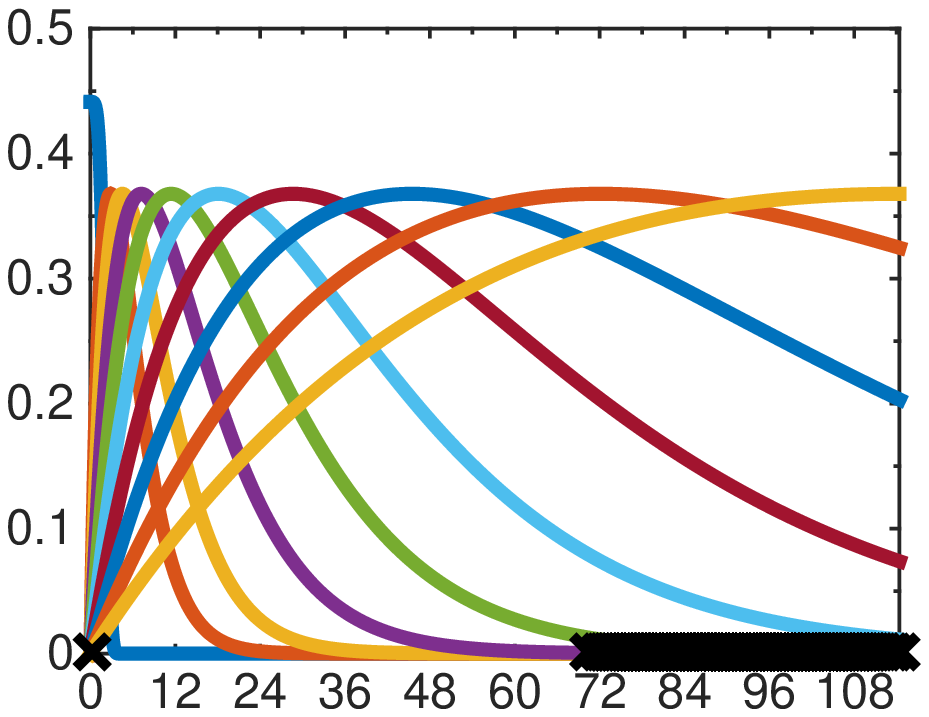} 
\includegraphics[width=0.45\textwidth]{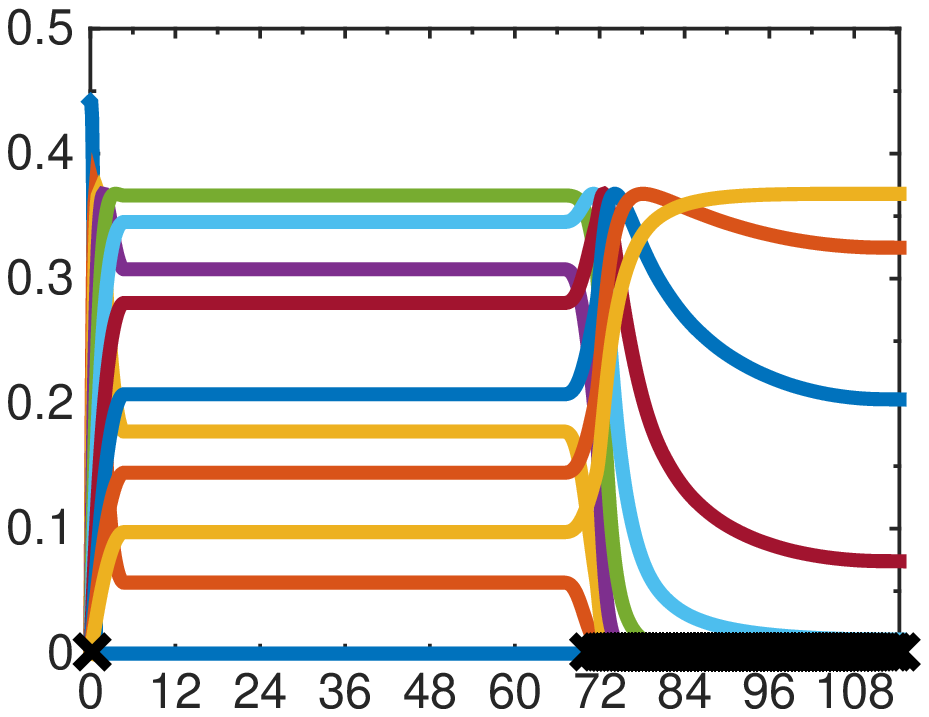} 
\end{center}
\caption{Non-adapted filterbank (left) and adapted filterbank (right).}
\vspace{-10pt}
\end{figure}
We notice that the Lanczos method exhibits excellent numerical behaviour in comparison with the Chebyshev method, particularly when using non-adapted filters. As the probability $p$ increases, the (relative) spectral gap increases as well, leading to more accurate Lanczos approximation (see Figure~\ref{fig:error_vs_p}).       



\section{Conclusion}
In this letter, a scalable, efficient and accurate method to perform signal filtering on graphs based on the Lanczos method is proposed. We compare it to the existing method based on the Chebyshev polynomial approximation. The advantage of filtering using the Lanczos approximation is that no a priori information about the spectrum of $\Lap$ is needed. In all numerical experiments the proposed method outperforms the Chebyshev approach. That is especially visible in case of a non-adapted filterbank. Futhermore, the Lanczos approximation is superior to the Chebyshev approximation when filtering on random graphs with well separated eigenvalues.

\bibliographystyle{abbrv}
\bibliography{biblio}

\end{document}